\documentstyle{amltd2004}
\begin{document}
\annalsline{156}{2002}
\received{January 20, 1999}
\startingpage{79}
\def\bye{\end{document}}
 \font\tenrm=cmr10
\def\dfrac#1#2{\displaystyle \frac{#1}{#2}}
\catcode`\@=11
\font\twelvemsb=msbm10 scaled 1100
\font\tenmsb=msbm10
\font\ninemsb=msbm10 scaled 800
\newfam\msbfam
\textfont\msbfam=\twelvemsb  \scriptfont\msbfam=\ninemsb
  \scriptscriptfont\msbfam=\ninemsb
\def\msb@{\hexnumber@\msbfam}
\def\Bbb{\relax\ifmmode\let\next\Bbb@\else
 \def\next{\errmessage{Use \string\Bbb\space only in math
mode}}\fi\next}
\def\Bbb@#1{{\Bbb@@{#1}}}
\def\Bbb@@#1{\fam\msbfam#1}
\catcode`\@=12

 \catcode`\@=11
\font\twelveeuf=eufm10 scaled 1100
\font\teneuf=eufm10
\font\nineeuf=eufm7 scaled 1100
\newfam\euffam
\textfont\euffam=\twelveeuf  \scriptfont\euffam=\teneuf
  \scriptscriptfont\euffam=\nineeuf
\def\euf@{\hexnumber@\euffam}
\def\frak{\relax\ifmmode\let\next\frak@\else
 \def\next{\errmessage{Use \string\frak\space only in math
mode}}\fi\next}
\def\frak@#1{{\frak@@{#1}}}
\def\frak@@#1{\fam\euffam#1}
\catcode`\@=12

\newcommand{\Prob}{{\rm Dist}}
\newcommand{\T}{{{\cal T}}}
\newcommand{\tid}{[T,{\rm Id}]}
\newcommand{\pendo}{uniform $p$\/{\rm -}\/to\/{\rm -}\/one endomorphism}
\newcommand{\pendos}{uniform $p$\/{\rm -}\/to\/{\rm -}\/one endomorphisms}
\newcommand{\X}{(X,T,\mu)}
\newcommand{\K}{\{0,\dots,p-1\}}
\newcommand{\Y}{(Y,S,\nu)}
\def\cont#1#2{\mathop{\subset}\limits^{#1}_{#2}}
\def\coup#1{\mathop{\vee}\limits^{#1}}

\title{Uniform endomorphisms which are\\ isomorphic to a Bernoulli shift}

\shorttitle{Uniform Bernoulli endomorphisms}  

 \twoauthors{Christopher Hoffman}{Daniel Rudolph}
 \institutions{University of Maryland, College Park, MD\\
{\eightpoint {\it E-mail address\/}: djr@math.umd.edu}\\
\vglue6pt
University of Washington, Seattle, WA\\
{\eightpoint {\it E-mail address\/}: hoffman@math.washington.edu}}

\centerline{\bf Abstract} 

A {\it uniformly $p$\/{\rm -}\/to\/{\rm -}\/one endomorphism} is a measure-preserving
map with entropy log $p$ which is almost everywhere $p$-to-one and for which the
conditional expectation of each preimage is precisely $1/p$.  The {\it
standard} example of this is  a one-sided $p$-shift with uniform i.i.d.\
Bernoulli measure.  We give a characterization of those uniformly
finite-to-one endomorphisms conjugate to this standard example by a
condition on the past tree of names which is analogous to {\it very weakly
Bernoulli} or {\it loosely Bernoulli.}  As a consequence we show that a large
class of isometric extensions of the standard example are conjugate to
it.
 
\section{Introduction}

Perhaps the most significant aspect of Ornstein's isomorphism theory
for Bernoulli shifts \cite{O1} 
is the very weak Bernoulli condition which Ornstein
and Weiss proved characterizes isomorphism to a Bernoulli shift \cite{O5}, \cite{OW2}.  
The very weak
Bernoulli (v.w.B.) condition has been exploited to show that many
classes of transformations are isomorphic to Bernoulli shifts.
Examples of these include ergodic toral automorphisms and geodesic
flows on negatively curved space \cite{K}, \cite{OW3}.  The Kakutani
equivalence theory of Feldman, Ornstein and Weiss has a similar
condition, loosely Bernoulli (l.B.), which shows when a transformation
is Kakutani equivalent to a Bernoulli shift \cite{ORW}.  The existence of such a
criterion is the hallmark of these and parallel theories.

We consider here {\it uniformly $p$\/{\rm -}\/to\/{\rm -}\/one endomorphisms}, 
measure-preserving endomorphisms with entropy $\log p$
which are a.s.\ $p$-to-one and for which the
conditional expectations of the preimages are all equal to $1/p$.  The
one-sided Bernoulli shift on $p$ symbols, each equally likely, is the
{\it standard} example of a uniformly $p$-to-one endomorphism.  
This endomorphism has state space $B=\{0,1, \dots,p-1\}^{{\Bbb N}}$.
The measure $\nu$ is defined by 
$$\nu(b \ | \ b_0 =a_0,b_1=a_1,\dots,b_m=a_m) =1/p^{m+1}$$
for any sequence  $a_0,\dots,a_m$ 
where all the $a_i \in \{0,1, \dots,p-1\}$.
The action on $B$ is the left shift $\sigma(b)_i=b_{i+1}$.
In this paper we give a characterization of those
uniformly $p$-to-one endomorphisms conjugate to the standard example
via a criterion parallel to v.w.B. and l.B.

We now describe the class of uniformly two-to-one
endomorphisms known as the 
$\tid$ examples.  Let $(B,\sigma,\nu)$ be the standard two-to-one endomorphism.
For $Y$ any Lebesgue
space and $T$ an ergodic automorphism of $Y$, $\tid$ is defined by
$$\tid (x,y)=(\sigma(x),T^{x_0}(y)).$$ 
Parry \cite{P} has shown that
for $T$, an irrational rotation of the circle extremely well
approximated by rationals, $\tid$ is always conjugate to the
standard example ($\sigma$ itself).  In Section \ref{tid} we will show 
that $\tid$ is conjugate to the standard example for all ergodic
isometries $T$.

The rest of the paper is organized as follows.
In Section \ref{twb} we define tree very weak Bernoulli and prove
that a large natural class of factors of tree v.w.B. endomorphisms are tree v.w.B. 
In Section \ref{osc} 
we will introduce the notion of a one-sided joining of two
endomorphisms. 
In Section \ref{tfd} we introduce the concept of tree finitely determined and
prove that tree very weak Bernoulli implies tree finitely determined.
Section \ref{cl} contains the bulk of the proof.
It will follow the same outline as the Burton-Rothstein approach to 
Ornstein's theorem.  First there will
 be a Rokhlin lemma and a strong Rokhlin lemma for endomorphisms.  This
 is followed by the copying lemma, which is the main tool in the proof.
 The proof of our copying lemma is much easier than
the proof of the copying lemma in  
Ornstein's theorem because, perhaps surprisingly, 
we will not have to deal with entropy.
In Section~\ref{tid} we apply the
tree v.w.B. characterization to show that a large class of isometric
extensions of the standard endomorphism are one-sidedly conjugate to it.

The methods of this work are much more broadly applicable than
to just uniformly finite-to-one endomorphisms.  For any bounded 
finite-to-one endomorphism $\X$ one can assign to each
point $x$ a tree of inverse images.  We will discuss this in
great detail as we continue.  Each point on this tree can
be weighted with the expectation that the path of inverse
images  to that point has occurred.  This assignment of
a weighted tree to $x$ is a conjugacy invariant
of the endomorphism.  Moreover, the expectations of just
the inverse images $T^{-1}(x)$ set a natural lower bound on the entropy
of $T$.  The uniformly finite-to-one maps are the simplest 
bounded  finite-to-one endomorphisms in that
the weighted trees do not depend on $x$ and the entropy of $T$ is
at its minimum given the tree.  To apply our methods to the general
bounded  finite-to-one case requires working relative to the 
sub-$\sigma$-algebra generated by the trees, with perhaps
nonzero entropy relative to the trees.  Both issues are standard
in the theory and require no essentially new ideas.  We have
focused on the simplest case here as it allows the presentation
of all the new ideas without the extra baggage and also as
it arises in a variety of natural situations, for example rational
maps of the Riemann sphere and ergodic group endomorphisms.
All of these issues will be presented separately. 

Isomorphism theories for endomorphisms have been presented
previously.  The two cases we are aware of are del Junco~\cite{J}
and Ashley, Marcus and Tuncel \cite{AMT}.  Del Junco considers two-sided
Bernoulli shifts and is after conjugacies that are finitary, one-sided in one
direction but two-sided in the other.  Thus his work is significantly
different from ours.  It is interesting to note his use of what
he calls a * joining which is strongly analogous to our one-sided
joinings.  The work of Ashley, Marcus and Tuncel is both
more restrictive and more general than ours as
they consider  all finite state Markov chains.  As our examples indicate,
many standard endomorphisms are not directly Markov.  On the
other hand, the existence of a period in a Markov chain will make
it nonstandard.   We expect future work on relativizing our methods
will provide an alternative approach to this classification of
one-sided Markov chains just as relativizing Ornstein's theorem
provided a classification of two-sided Markov chains.

 The authors wish to thank A. del Junco for discovering several gaps
and one serious error in the previous version of this paper.

\section{Tree very weak Bernoulli}
\label{twb}

In this section we will define what it means for an endomorphism to be
tree very weakly Bernoulli.
As with v.w.B. and l.B., our criterion concerns names on orbits.   
A $p$-to-$1$ uniform endomorphism $T$ and a point $x$ generate a
$p$-ary tree of inverse images.  We will describe how to define 
a {\it tree name} given the tree of inverse images and 
any partition or, more generally,
any metric space-valued function on $X$.  
Simply stated the 
tree very weakly Bernoulli
criterion says that for any such function
and a.e.\ two points $x$ and $x'$ one can match these two tree names
with an arbitrarily small density of errors by a map that preserves the
tree structure.  

Consider a $p$-ary tree with $p^{n}$ nodes at each index $n\geq 0$.  Each
node at index $n$ connects to  $p$ nodes at the index $n+1$.  We assign each
such set of $p$ nodes a distinct value in $\{0,\ldots ,p-1\}$.  Then we
label each node other than the root 
by the sequence of values we see moving from the root
to the node.
In this form we can concatenate nodes $v'$ and $v$
by concatenating their labels.  Notice that when we fix a node $v$, the
set of labels $vv'$ form a $p$-ary subtree rooted at node $v$.  This
is consistent with our convention that the root node is unlabeled.

Call this labeled tree ${{\cal T}}$ and call the
tree that has all the nodes at index less than or equal to $n$,
${{\cal T}}_n$.
Let $\eta$ be the set of nodes of ${{\cal T}}$ and $\eta_n$ be the set of nodes for
${{\cal T}}_n$.  
For $v\in \eta$ and at index $i$ 
(i.e.\ $v \in \eta_i \setminus \eta_{i-1}$) we write $|v|=i$ and where $v$ is a list of values
$v_1,\dots,v_{i}$ from $\{0,\ldots ,p-1\}$ 
which is the list of labels of the nodes
along the branch from the root to $v$.

Let ${\cal A}$ be the collection of all bijections of the nodes of
${\cal T}$ that preserve the tree structure.
We refer to this as the group of {\it tree automorphisms}.  Let
${{\cal A}}_n$ be the bijections of the nodes of ${{\cal T}}_n$
preserving the tree structure.  To give a representation to such
automorphisms $A$ notice that from $A$ we obtain a permutation $\pi_v$
of $\{0,\ldots ,p-1\}$ at each node, giving the rearrangement of its 
$p$ predecessors.
An automorphism
of ${{\cal T}}_n$ will be represented by an assignment of a
permutation of $\K$ to each node of the tree 
except for those at index $n$.

Let $(X,T,\mu, {\cal F})$ be a uniformly $p$-to-one endomorphism.
Then each $x\in X$ has $p$ inverse images.
Select a measurable $p$ set partition $K$ of $X$
such that almost every $x$ has one preimage in each element of $K$.
Label the sets of $K$ as $K_0,K_1,\ldots ,K_{p-1}$. 
For each $i\in \{0,\ldots ,p-1\}$ and $x \in X$ define $T_i(x)$ to
be the preimage of $x$ in $K_i$.  
We now define a set of {\it partial inverses} for $T$.
For $v=(v_1,\dots ,v_{i}) \in  \eta$
define $T_v(x)=T_{v_{i}}(\cdots(T_{v_{1}}(x)))$.

 We let $R$ and $U$ be compact metric spaces.
We will  use $d$ generically for a metric on $R$ and $U$ and
more generally for any metric space  considered,
and will always assume that the labeling spaces $R$ and $U$ have $d$ diameter precisely~$1$.

A ${{\cal T}},R$ name $h$ is any function from ${{\cal T}}$ to $R$.
We say it is {\it tree adapted} if for any $v \in \eta$ and $i,j \in \K$
with $i \neq j$ we have $h(vi)\neq h(vj)$ and we  say it is {\it strongly tree adapted}
if $d(h(vi), h(vj))=1$.  A map $f: X \to R$ generates
${{\cal T}},R$ names by ${{\cal T}}_x(v)=f(T_v(x)).$
We say that $f$ is {\it {\rm (}\/strongly\/{\rm )} tree adapted} if 
${{\cal T}}_x(v)$ is (strongly) tree adapted for almost every $x$.  
Let ${\cal G}$ be the\break $\sigma$-algebra on $X$ generated by pullback of the
$\sigma$-algebra on $R$.  We say $f$
generates if $\vee _i T^{-i}({\cal G}) = {\cal F}$. 
``Strongly tree adapted" depends explicitly on the choice of metric but 
obviously, tree adapted, does not. 
For any $f:X\to R$, the map 
$f\vee K:X\to R\times\{0,\dots,p-1\}$ will always
be tree strongly  tree adapted.  More generally so long as $d(h(vi),h(vj))$ is
bounded uniformly away from zero by some $\alpha>0$ we can replace $d$ by an equivalent
metric of diameter 1 making $h$ strongly tree adapted.  Just set
$$d_{\rm new}(x,y)=\min(1, d(x,y)/\beta)$$ for some $0<\beta\leq \alpha$.

Any measure-preserving endomorphism $(X,T,\mu)$ and function $f:X\to R$ 
generate  a stationary sequence of random variables $R_i=f\circ
T^i, i\geq 0$.  We regard such a sequence as \pagebreak a measure on $R^{{\Bbb N}}$
with the weak* topology.  There is a unique extension of this measure to all of
$R^{{\Bbb Z}}$, preserving stationarity.  For a.e. $x\in X$ let
$\{R_i(x)\}_{i<0}$ be random variables with distribution
${\rm Dist}(R_i, i<0|R_i=f(T^i(x)), i\geq 0)$.

Any tree adapted ${{\cal T}}, R$-name $h$, generates an $R$-valued
sequence of random variables $\{R_i(h)\}_{i<0}$. To a cylinder set
$r_{-j},r_{-j+1},\dots,r_{-1}$ we assign the measure equal to $p^{-j}$
times the number of all nodes at index $j$ whose name to index $1$ is
the word $r_{-j},\dots,r_{-1}$.  For $(X,T,\mu)$ a uniform $p$-to-one
endomorphism and $f:X\to R$ tree adapted and generating and
${{\cal T}}_x$ the ${{\cal T}},R$-name of $x$, $\{R_i(x)\}_{i<0}$
and $\{R_i({{\cal T}}_x)\}_{i<0}$ are just two descriptions of the
same sequence of random variables.

We now put a family of metrics on ${{\cal T}},R$-names (and on ${{\cal T}}_n,R$-names).  
For two ${{\cal T}},R$-names $h$ and $h'$ we define
$$\overline t_n(h,h')=\min_{A\in{{\cal A}}_n}\frac1{n} \sum_{0<|v|\leq n}
	p^{-|v|}d(h(v),h'(A(v))).$$

For the two names $h$ and $h'$ in this definition, if one
follows each branch from the root and through the tree and writes down
the name seen along that branch, one obtains $p^{n}$ different names.
Giving each name a mass of $p^{-n}$ one obtains the two
sequences of random variables $R_i(h)$ and $R_i(h'), -n\leq i <0$. 
 A matching of the trees via a
tree automorphism gives a coupling of these two distributions.  The weighting
of nodes is such that the calculation of the $\overline t$ distance is
precisely the $\overline d$ distance between the  random variables
one would calculate from this
coupling.  In this sense $\overline t$ is at least as large as
$\overline d$, which would be the inf over all couplings, not just
those that come from tree automorphisms.

\numbereddemo{Definition} 
Let $\X$ be a \pendo\ and
$f$ a tree adapted map from $X \to R$. 
We say $\X$ and $f$ are {\it tree very weak Bernoulli {\rm (}\/tree v.w.B.}\/) if
 for
any $\varepsilon>0$ and all $N$ sufficiently large there is a set $G=G(\varepsilon, N)$ 
with
$\mu(G)> 1 - \varepsilon$ such that for any $x,y \in G$,
$$\overline t_N({{\cal T}}_x,{{\cal T}}_y)< \varepsilon.$$
\enddemo
 
It is fairly direct that tree v.w.B. endomorphisms are exact and hence always
ergodic \cite{H}, \cite{V}. 
Our next goal is  to show that a large natural class of factors of tree v.w.B. endomorphisms is again tree v.w.B. First we
define the class we are interested in.

\numbereddemo{Definition} 
We say a factor map $\phi$ from $(X,T,\mu)$ to 
$(Y,S,\nu)$ is
{\it tree adapted} if for a.e.\ point $x$ the map $\phi$ restricted
to $T^{-1}(x)$ is one-to-one {\it into} 
$S^{-1}(\phi(x))$. 
\enddemo 

It is not difficult to construct factor maps between endomorphisms that are
not tree adapted.  On the other hand a 
conjugacy is clearly tree adapted in both directions.  Our
constructions will never leave the class of tree adapted factor maps.
We now give a useful little technical lemma that tree
adapted factors of uniform $p$-to-one endomorphisms are themselves
uniform $p$-to-one endomorphisms.

\proclaim{Lemma} \label{pendotopendo}
Suppose $\X$ is a \pendo\ and $\Y$ is a
tree adapted factor of $\X$ by a map $\phi$.  Then $\phi$ restricted
to $T^{-1}(x)$ is almost surely onto $S^{-1}(\phi(x))$ and $\Y$ is also a \pendo.
\endproclaim 

\demo{Proof} 
Because the map $\phi$ is tree adapted the endomorphism $Y$ is a.e.\
at least $p$-to-one.  For each $x \in \phi^{-1}(y)$, $x'\in T^{-1}(x)$ and $\phi(x')=y'$
  the conditional probability 
of $x'$ given $x$ is the same as the conditional probability of $y'$
given $x$ and $y$ (i.e.\ $1/p$).    Thus
the entropy is at least log $p$.  As the entropy of $X$ is log $p$ the
entropy of $Y$ can be at most log $p$.  Thus the conditional probability
of $y'$ given only $y$ must be $1/p$ and $\Y$ is a \pendo. 
\enddemo 

\proclaim{{C}orollary} 
  For $(X, T,\mu)$ a \pendo{\rm ,} $R$ compact metric and $f:X\to R$ tree adapted{\rm ,}
the left shift on the sequence of random variables $R_i=f\circ T^i$ is a \pendo.
\endproclaim

\proclaim{Lemma} 
\label{factor} \hskip-8pt
Suppose $\X$ is a \pendo{\rm ,} $f$~is tree adapted and generating{\rm ,} 
$\X$ and $f$ are tree {\rm v.w.B}.\ and  $\Y$ is a  
factor by a tree adapted factor map $\phi$.
Then for any tree adapted 
$g:Y \to R${\rm ,} $\Y$ and $g$ are also tree {\rm v.w.B}.
\endproclaim

 \demo{Proof} 
By Lemma 2.3, we know $\Y$ is a \pendo.
Let $\cal G$ be the $\sigma$-algebra of $f$ measurable
sets; hence $g\circ\phi$ is $\bigvee_{i=0}^\infty T^{-i}(\cal G)$
measurable.  Writing this as a finite approximation, for each $\varepsilon>0$
we can find an $s\in{\Bbb N}$, a $\delta>0$ and a subset $G'\subseteq X$
with $\mu(G)>1-\varepsilon$ so that if $x,x'\in G'$ and
$$d(f(T^i(x)),f(T^i(x')))<\delta\hbox{ for }0\leq i<s\hbox{ then}$$
$$d(g(\phi(x)),g(\phi(x')))<\varepsilon.$$

As $T$ is ergodic the mean ergodic theorem tells us that
$$\frac 1n\sum_{0<|v|\leq n}p^{-|v|}\chi_{G'}(T_v(x))\stackrel{n}{\to}\mu(G')$$
in $L^1$. (Convergence here can be shown to be pointwise but as this
is not standard and we do not need it we just quote the mean convergence which
follows directly from rewriting this average as an average over the forward images
of $T^{-n}(x)$.)

Suppose $x$, $x'$ and $N$ satisfy
$$\overline
t_N({{\cal T}}_x,{{\cal T}}_{x'})<\dfrac{\delta^2}{2s},
$$
$$\frac 1N\sum_{0<|v|\leq N} p^{-|v|}\chi_{G'}(T_v(x))>1-2\varepsilon
$$ and
$$\frac 1N\sum_{0<|v|\leq N} p^{-|v|}\chi_{G'}(T_v(x'))>1-2\varepsilon.$$
Then there is an $A\in{{\cal A}}_n$ given by
$$\frac 1N\sum_{0<|v|\leq N}p^{-|v|}d(f(T_v(x)),f(T_{A(v)}(x')))<\frac{\delta^2}{2s}.$$

Let $Z=\{v|0<|v|\leq N, d(f(T^i\circ T_v(x)),f(T^i\circ T_{A(v)}(x')))<\delta\}$ and 
conclude
$$\frac 1N\sum_{v\in Z}p^{-|v|}>1-\delta.$$
If $v\in Z$ and both $T_v(x)$ and $T_{A(v)}(x')$ are in $G'$ then 
$$d(g\circ\phi(T_v(x)), g\circ\phi(T_{A(v)}(x')))<\varepsilon
	\hbox{ and hence}$$
$$\frac 1N\sum_{0<|v|\leq N} 
p^{-|v|}d(g(S_v(\phi(x))),g(S_{A(v)}(\phi(x'))))<5\varepsilon+\delta.$$
That $(Y,S,\nu)$ is tree v.w.B. now follows.
\enddemo 

 \vglue-12pt

 \section{One-sided couplings}\label{osc}
\vglue-6pt

As we are following the Burton-Rothstein approach to the isomorphism
theorem we will be considering
joinings of two endomorphisms $(X,T,\mu)$ and $(Y,S,\nu)$.
A {\it coupling} of two spaces $(X,\mu)$ and $(Y,\nu)$ is a 
measure on 
$X \times Y$ which has marginals $\mu$ and $\nu$.  A {\it joining} of 
$(X,T,\mu)$ and $(Y,S,\mu)$ is a coupling of $(X,\mu)$ and
$(Y,\nu)$ which is invariant under $T \times S$.
We will not consider all of the joinings of $(X,\mu)$ and $(Y,\nu)$
but rather a collection we call one-sided joinings.  In this section we
define one-sided joinings and prove some facts about them.

\numbereddemo{Definition} 
Suppose $S$ is a subset of ${\Bbb Z}$ and
$\{R_i\}_{S}$ and $\{U_i\}_{S}$ 
are two sequences of random variables.
  A {\it one\/{\rm -}\/sided coupling} is a
coupling of these two sequences such that for all $n$ and $i>n$
\begin{eqnarray*}
\Prob(\{R_j\}_{j< i}|\{R_j\}_{ j\geq i},\{U_j\}_{ j\geq i})&=&\Prob(\{R_j\}_{j< i}|\{R_j\}_{j\geq 
i}),\\
\noalign{\noindent and symmetrically}
\noalign{\vskip6pt}
  \Prob(\{U_j\}_{j< i}|\{R_j\}_{j\geq i},\{U_j\}_{j\geq i})&=&\Prob(\{U_j\}_{j< 
i}|\{U_j\}_{j\geq i}).
\end{eqnarray*}
If the sequences are stationary 
then the one-sided couplings that also are stationary are called
{\it one\/{\rm -}\/sided joinings.} 
\enddemo 

If along with the endomorphism $(X,T,\mu)$ there is a 
generating function
$f:X \to R$  we have defined $R_i(x)=f(T^i(x))$ for any $x$ and $i \geq 0$.
We extend this by stationarity to negative $i$ giving a stationary
sequence $R_i$, $i\in {\Bbb Z}$.  (Note that the map from stationary
measures on $R^{{\Bbb N}}$ to those on $R^{-{\Bbb N}}$ obtained
by extension to ${\Bbb Z}$ and then restriction is a weak* homeomorphism.)
If we also have an endomorphism $(Y,S,\nu)$ and a 
generating function
$g:Y \to R$  
we define $U_i(y)=g(S^i(y))$ and extend to indices $i<0$ similarly.

\numbereddemo{Definition} 
A coupling (or a joining) of $\X$ and $\Y$ is {\it one\/{\rm -}\/sided} if there are 
generating functions $f$ and $g$ so that the coupling of 
$R_i$ and $U_i$, $i\geq 0$,  is one-sided.
\enddemo 

Notice that if a coupling   generates a one-sided  coupling
for one choice of generating functions $f$
and $g$ it will do so for all choices. 

\proclaim{Lemma}   \label{closedcond}
  For sequences of random variables $R_j$ and $U_j$ and each value~$i${\rm ,}
those couplings for which
$$\Prob(\{R_j\}_{j< i}|\{R_j\}_{ j\geq i},\{U_j\}_{ j\geq i})=\Prob(\{R_j\}_{j< i}|\{R_j\}_{j\geq 
i})$$
are a weak* closed subset of all couplings.
\endproclaim

\demo{Proof}  We begin with some basic reductions of the problem.
  Notice that this statement is simply  about couplings of three
measure algebras $X_1=\vee_{j<i}R_j$, $X_2=\vee_{j\geq i}R_j$ and $Y=\vee_{j\geq i}U_i$
where the first two are coupled by a fixed measure $\mu$.  We consider those couplings
of the three where $X_1$ and $Y$ are coupled conditionally, independently over $X_2$.
Viewed this way we see that without loss of generality we can assume all the random variables $R_i$ and
$U_i$  (to be explicit) are two-valued.  Next notice that the condition reduces to the countable
list of conditions

 $$\Prob(\{R_j\}_{I\leq j< i}|\{R_j\}_{ j\geq i},\{U_j\}_{ j\geq i})=\Prob(\{R_j\}_{I\leq j< i}|\{R_j\}_{j\geq 
i}).$$
To show that each such condition is weak* closed is equivalent to proving closedness when $X_1$ is a finite space,
i.e.\ we have reduced to the case of showing closedness of
$\Prob(P|X_2,Y)=\Prob(P|X_2)$ where $P$ is a finite partition.  Now suppose $\nu_i$ all satisfy this
condition and converge weak* to $\nu$.  Upper semi-continuity of entropy tells us that
$\limsup H_{\nu_i}(P|X_2\vee Y)\leq H_{\nu}(P|X_2,Y)$.  On the other hand
$$H_{\nu_i}(P|X_2\vee Y)=H_{\mu}(P|X_2)= H_{\nu}(P|X_2)\geq H_{\nu}(P|X_2\vee Y)$$
and $H_{\nu}(P|X_2\vee Y)=H_{\mu}(P|X_2)$ and $\mu$ must also satisfy the condition. 
\enddemo 
\pagebreak

\proclaim{{C}orollary} \label{partialsgroup} Suppose $\mu_1$ couples $U_i$ and $R_i$
 one\/{\rm -}\/sidedly and $\mu_2$ couples
  $V_i$ and $R_i$ one\/{\rm -}\/sidedly.  There is then a coupling $\hat\mu$ of all three sequences
that projects to $\mu_i$ on the appropriate pair of sequences and when restricted to
$U_i$ and $V_i$ is one\/{\rm -}\/sided.
\endproclaim

\demo{Proof}   Suppose $\overline\mu$ couples the three sequences $R_i$, $U_i$ and $V_i$
  for $i\geq I$ with the given marginals in some way, not necessarily one-sidedly.  Extend
$\overline \mu$ to index $I-1$ by first coupling on $R_{I-1}$ relatively independently over
the algebra $\bigvee _{i\geq I}R_i$.  By the one-sidedness of the marginals, this maintains the
marginals.  Now couple on $U_{I-1}$ to this relatively independently over the algebra
$\bigvee_{i\geq I-1}R_i \vee \bigvee_{i\geq I}U_i$.  Couple on $V_{I-1}$ symmetrically.
This preserves the marginals and  
$$\Prob(U_{I-1}|\{U_i\}_{i\geq I},\{V_i\}_{i\geq I})=\Prob(U_{I-1}|\{U_i\}_{i\geq I}\}),$$
and symmetrically,
$$\Prob(V_{I-1}|\{U_i\}_{i\geq I},\{V_i\}_{i\geq I})=\Prob(V_{I-1}|\{V_i\}_{i\geq I}\});$$
i.e., for this one step the coupling is one-sided.

We continue inductively to add on variables as $I\to -\infty$ and   obtain a coupling which
is one-sided at all indices $i\leq I$.  For each value $I\geq 0$ start with $\overline\mu$,
the relatively independent coupling of the $\mu_i$ and $\mu_2$ over the common
algebra $\bigvee_{i\geq I} R_i$.  Extend to the right as described above to produce
a measure $\overline\mu_I$.  Let $\hat\mu$ be any weak* limit of the couplings
$\overline\mu_I$.  Lemma~\ref{closedcond} guarantees that $\hat\mu$ is one-sided.  
\enddemo

\numbereddemo{Definition} 
 Suppose $(X,T,\mu)$ and $(Y,S,\nu)$ are two endomorphisms. 
 Then   define $J^+((X,T,\mu),(Y,S,\nu))$ to be the space
 of all one\/{\rm -}\/sided joinings of these endomorphisms.
 Define $J^+_e((X,T,\mu),(Y,S,\nu))$ to be the {\it ergodic} and
 one\/{\rm -}\/sided joinings.  If there is no confusion we will just write
 $J^+$ and $J_e^+$.
\enddemo 

Notice that Lemma~\ref{closedcond} has shown the one-sided joinings to
be a closed subset of all joinings.   It is not difficult to see that if $T$ and
$S$ are assumed ergodic, then the one-sided joinings are convex and
that the extreme points of this set are the ergodic and one-sided joinings.
(One just notes that the ergodic components of a one-sided joining must
themselves be one-sided as the ergodic decomposition is past measurable.)
Furthermore we see that if $\hat\mu_1\in J^+$
of $T$ and $S$ and $\hat\mu_2\in J^+$ of $T$ and $U$ then there is
a stationary joining of all three which projects to these two on the
appropriate pairs and is in $J^+$ of $S$ and $U$ on this pair.  Just
observe that Corollary~\ref{partialsgroup} gives us a one-sided coupling
and by averaging over translates and taking a weak* limit we get a joining.
By extending this further, we see that  if the $\hat\mu_i$ were in $J^+_e$, i.e.\ were ergodic,
then the one-sided joining of $S$ and $U$ can also be chosen ergodic.
Almost any ergodic component of the one-sided joining just constructed will do.

We consider two different weak* pseudometrics on processes of the form
$(\X,f)$.
First an endomorphism $\X$ and a function $f:X \to R$ define a measure on 
$R^{{\Bbb N}}$.   Let ${\rm dist}$ be a metric 
on $C^*(R^{{\Bbb N}})$ which generates the weak* topology.
When we refer to 
$${\rm dist}((\X,f),(\Y,g))={\rm dist}(f,g)$$ 
we mean the   
distance between the measures that $(\X,f)$ and $(\Y,g)$ generate on 
$R^{{\Bbb N}}$.
We mention a particular case of this ${\rm dist}$  (pseudo)metric topology  to be
used repeatedly.
If $\hat \mu$ is a joining of $\X$ and $\Y$ and $f:X\to R$ and 
$g:Y\to U$ then $\hat \mu$ projects to a stationary measure on 
$(R \times U)^{{\Bbb N}}$.  
In this case we use the notation
$${\rm dist}(f\coup{\hat \mu} g,{\tilde f}\coup{\nu} {\tilde g}). $$

A \pendo \  $\X$ and a function $f:X \to R$ define  a measure on 
$R^{{\cal T}}$ as we have associated with each point $x
\in X$ the ${{\cal T}},R$ name  ${{\cal T}}_x$.  
Let ${\rm tdist}$ be a metric 
 for the weak* topology on Borel measures on $R^{{{\cal T}}}$.
When we refer to
$${\rm tdist}((\X,f),(\Y,g))={\rm tdist}(f,g)$$ we mean the ${\rm tdist}$
distance between the measures these processes generate on~$R^{{\cal T}}$.

\proclaim{Lemma} \label{dististdist}
  The {\rm (}\/pseudo\/{\rm )}\/topologies generated by ${\rm dist}$ and ${\rm tdist}$
on uniform $p$\/{\rm -}\/to\/{\rm -}\/one endomorphisms and tree adapted 
functions to $R$ are the same.
\endproclaim

 \demo{Proof} By Lemma~\ref{pendotopendo} we can assume $f$ generates.
We have seen any tree adapted ${{\cal T}}, R$ name $h$ generates a sequence
of random variables $R_i(h)$, i.e.\ a measure on $R^{-{\Bbb N}}$.
 Thus any measure $m$ on the tree adapted ${{\cal T}}, R$ names
projects to a measure on $R^{-{\Bbb N}}$.  
If this measure comes from a \pendo \ then it is stationary and
maps homeomorphically to a measure on $R^{{\Bbb N}}$.
Call this measure $\Psi(m)$. Obviously, $\Psi$ is
 a weak* continuous map which shows that ${\rm tdist}$ is at
least as strong as ${\rm dist}$.

For a \pendo\ $\X$ and  a.e.\ $x$, the map\break
$x\to {{\cal T}}_x$ 
lifted to measures is an inverse for $\Psi$.  

To see  that $\Psi^{-1}$ is continuous, suppose $\{R_i\}$ and $\{R_i^j\}$, $j=1,2,\ldots$,
are sequences of random variables (measures on $R^{{\Bbb Z}}$) arising from
tree-adapted functions on uniformly $p$-to-one endomorphisms with
$\{R_i^j\}\longrightarrow \{R_i\}$ in dist.  This is equivalent to putting all these
random variables on a common measure space $(\Omega,\hat\mu)$ (which might as well
be $(R^{{\Bbb N}})^{{\Bbb Z}}\times R^{{\Bbb Z}}$) with each
$R_i^j\longrightarrow R_i$ in $\hat\mu$ probability.\pagebreak

We write $\vec R=\{R_i\}_{i\in {\Bbb Z}}$, $\vec R^j=\{R_i^j\}_{i\in {\Bbb Z}}$ and
$\vec R_+=\{R_i\}_{i\in{\Bbb N}}$ etc.  Let ${\cal T}_{\vec R_+}$ be the uniform $p$-adic tree of inverse
images of $\vec R_+$ and ${\cal T}_{\vec R^j_+}$ that of $\vec R_+^j$.  All these labeled trees have the property
that the labels of the $p$ predecessors of any node are precisely 1 apart in $R$.  The collection of such $R$ labelings
of the tree are a closed subset $C$ of $R^{{\cal T}}$.

Let $C_n\subseteq R^{{\cal T}_n}$ consist of the labelings in $C$ restricted to ${\cal T}_n$.  On $R^{{\cal T}_n}$ use
the sup metric up to tree automorphisms and on $R^n$ the sup metric.

A labeling $\eta\in C_n$ gives rise to $p^n$ distinct names in $R^n$ -- the names along the $p^n$ branches.
Call this set of names $N(\eta)$.  The critical observation here is this:  Suppose for two labelings $\eta_1$ and $\eta_2\in C_n$,
each element of $N(\eta_1)$ is within $\varepsilon <1/2$ of some element in $N(\eta_2)$.  Then the labelings $\eta_1$
and $\eta_2$ themselves must be within $\varepsilon$.  Just notice that the labels along distinct branches of $\eta_1$ cannot
be matched within $\varepsilon$ of the same branch of $\eta_2$ and moreover the matching must preserve the tree structure.

Both $R^{{\Bbb N}}$ and $R^{{\cal T}_n}$ are compact and so for each $\varepsilon>0$ there is a closed subset $G=G(\varepsilon)$
with $\hat\mu(G)>1-\varepsilon^2/p^n$ so that the map $\vec R_+\to {\cal T}_{\vec R_+,n}$ is uniformly continuous on $G$.  Hence there
are an $N$ and $\delta$ so that if $\vec R_+$ and $\vec R'_+\in G$ and
$d(R_i, R'_i)<\delta$ for $i\leq N$ then $d({\cal T}_{\vec R_+,n},{\cal T}_{\vec R'_+,n})<\varepsilon/2.$

As $\vec R^j\longrightarrow \vec R$ in probability, for all $j$ large enough there will be a subset $H=H(\varepsilon,j)$ of values
$\vec R^j$ with $\mu(H)>1-2\varepsilon/p^n$ and for each $\vec R^j\in H$ there is a {\it representative value}
$\vec R(\vec R^j)$ with $\vec R_+(\vec R^j)\in G$   so that

\begin{itemize}
\item[i)] for $-n\leq i\leq N$ we have $d(R_i^j,R_i(\vec R^j)<\varepsilon/2$ and
\vspace{.1in}
\item[ii)] $E(\vec R_+\in G\quad{\rm  and }\quad d(R_i,R_i(\vec R^j))<\delta,\quad 1\leq i\leq n|\vec R^j)>1-\varepsilon.$
\end{itemize}

Let $H'\subseteq H$ consist of those $\vec R^j$ for which all $p^n$ extensions of $\vec R_+^j$ in 
$N({\cal T}^j_{\vec R_+^j})\circ \vec R_+^j$
intersect $H$ nontrivially, i.e.\ have {\it good representatives}.  Now $\mu(H')>1-2\varepsilon$ and for any
$\vec R_+^j\in H'$ all names in $N({\cal T}_{\vec R_+^j,n}^j)$ must be within $\varepsilon$ of some branch of
${\cal T}_{\vec R_+(\vec R^j),n}$ and by our observation above
$$d\bigl({\cal T}_{\vec R_+^j,n},{\cal T}_{\vec R_+(\vec R^j),n}\bigr)<\varepsilon.$$ 
Now, by ii), for $\vec R^j\in H'$ 
we have
$$E(d({\cal T}_{\vec R^j_+},{\cal T}_{\vec R_+})<2\varepsilon|\vec R^j)>1-\varepsilon$$
and so
$$E(d({\cal T}_{\vec R_+^j,n},{\cal T}_{\vec R_+,n})<2\varepsilon)>1-3\varepsilon.$$
We conclude that $\vec R^j\to \vec R$ in tdist.
\enddemo

This next lemma is important because it says that all of the joinings
  created in the next section are one-sided.

\proclaim{Lemma} 
Suppose $(X,T,\mu)$ is a uniform $p$\/{\rm -}\/to\/{\rm -}\/one endomorphism 
and the factor map
  $\phi$ to $(Y,S,\nu)$ is  tree adapted.  Then the joining
$\hat\mu$ 
of $(X,T,\mu)$ and $(Y,S,\nu)$ generated by $\phi$ is one\/{\rm -}\/sided.
\endproclaim

 \demo{Proof} 
Fix $x$ and $y=\phi(x)$.
By stationarity it is sufficient to show that the conditional
probability of preimages of $y$ is the same as the conditional
probability of preimages of $y$ given $x$ and $y$.
As $\phi$ is a bijection from the inverse
images of $x$ to those of $y$ the conditional mass of any $T_v(x)$
given $x$ must
be precisely that of $S_v(y)$ given $y$.  This value
is $p^{-|v|}$ and hence the conditional expectation of each $S_{v'}(y)$
given 
$x$ and $y$ is the same as its expectation given $y$.  The other
set of equalities  is obvious as $x$ determines $y$ so  that
conditioning on $x$ and $y$ is the same as 
conditioning on $x$.
\enddemo

\numbereddemo{Definition} 
Suppose $(X,T,\mu)$ and $(Y,S,\nu)$ are two \pendos\ and $f$ and $g$
are functions into the same metric space $(R,d)$ and $n>0$.  Now,
\begin{eqnarray*}
\hat t_n(f,g)&=&\hat t_n((\X,f),(\Y,g))\\
&=&
	\inf\left(\frac 1n\sum_{i=1}^{n} 
    \int d(f(T^i(x)), g(S^i(y))\,d\hat \mu \right)\end{eqnarray*}
 and
\begin{eqnarray*}
\bar t_n(f,g)&=&\bar t_n((\X,f),(\Y,g))\\
&=&\inf\left( 
    \int \bar t_n({{\cal T}}_x,{{\cal T}}_y)\,d\hat\mu\right),\end{eqnarray*}
where the {\it inf}\/'s are taken over all one-sided couplings $\hat\mu$.
For comparison's sake we include the definition
\begin{eqnarray*}
\overline d_n(f,g)&=&\overline d_n((\X,f),(\Y,g))\\
&=&
	\inf\left(\frac 1n\sum_{i=1}^{n}
    \int d(f(T^i(x)), g(S^i(y))\,d\hat\mu \right),
\end{eqnarray*}
where the {\it inf} is taken over all couplings $\hat\mu$.
\enddemo 

On the face of it, $\bar t_n$ is a metric on random variables indexed on $-n\leq i<0$ and
$\hat t_n$ on random variables indexed on $0\leq i<n$.  By stationarity of $R_i$ these
can be translated to be the same sets of random variables.
  The proof of the following lemma will be given later in the
section.

\proclaim{Lemma} \label{t-hatequalst-bar}
  Suppose $\X$ and $\Y${\rm ,} are \pendos\ with tree adapted functions 
$f$ and $g$ to $R$.  Then
$\bar t_n(f,g)=\hat t_n(f,g)$.
\endproclaim 

We will not use the notation $\hat t_n$ again except in the proof
of Lemma~\ref{t-hatequalst-bar}, using just $\overline t_n$ for both
notions.

 \numbereddemo{Definition} 
   Suppose $\X$ and $\Y$, are 
\pendos\ and $f$ and $g$ are functions with values in 
$R$.  We set
 $$\overline t(f,g)=\liminf_{n\to\infty}\overline t_n(f,g).$$
 \enddemo  

\proclaim{Lemma}  
   Suppose $\X$ and $\Y$ are two \pendos\ and $f$ and $g$ are functions
   to $R$.  There is then a stationary{\rm ,} ergodic
   and one\/{\rm -}\/sided joining{\rm ,} $\mu,$ with $$\overline t(f,g)=\int
   d(f(x),g(y))\,d\mu.$$ In particular the $\liminf$ in the
   definition is a limit.  
\endproclaim

 \demo{Proof}  
   The simple weak* compactness argument completely analogous to that for
   $\overline d$  works here as the set of one-sided couplings
   is closed and convex, and the extreme points of the stationary and
   one-sided couplings are the ergodic ones.  \enddemo 

\numbereddemo{Definition} 
  Suppose $(X,T,\mu)$ and $(Y,S,\nu)$ 
are two endomorphisms,
$f:X\to R$ is tree-adapted and $\hat \mu\in J^+_e$ is an ergodic one-sided joining
of them.  We
say $f\cont{\varepsilon}{\hat\mu}Y$ if there is a one-sided and tree-adapted
function $\bar f :Y \to R$ with
$$\int d(f(x),\bar f(y))\,d{\hat \mu}{(x,y)} < \varepsilon.$$
Now, $f\cont{0}{\hat\mu}Y$ if $f\cont{\varepsilon}{\hat\mu}Y$ for
all $\varepsilon >0$.
\enddemo
 
Notice that if $f\cont{\varepsilon}{\hat\mu} Y$ one immediately obtains that for  
$\bar f$ of the definition  
$$\bar t((\X,f),(\Y,\bar f))<\varepsilon.$$
\noindent
Also notice that if $f$ generates and 
$f\cont{0}{\hat\mu}Y$ then relative to $\mu$ the endomorphism $\X$
sits as a one-sided and tree-adapted factor of $\Y$.

Now we show that a one-sided coupling lifts naturally to a measure on 
$X \times Y \times {\cal A}$.  This is the essential ingredient for showing that
$\hat t$ and $\bar t$ are equal and is also necessary in 
the copying lemma.
This lift is   not unique.   
(In the form we now describe, the direct product
of two \pendos\ has many potential lifts to a third automorphism coordinate.)
Let $h$ and $h'$ be
 $R$-valued and $U$-valued tree names and $R_i(h)$ and $U_i(h'), i<0$, be the
sequences of random variables they generate.  \pagebreak For any automorphism 
$A$ we construct a joined name
$\hat h_A(v)=(h(A^{-1}(v)),h'(v))$.  
Such a name will project to a measure on $R^{-{\Bbb N}}\times U^{-{\Bbb N}}$
that is a one-sided coupling of $R_i(h)$ and $U_i(h')$.   
Call it $\hat \mu_{(h,h',A)}.$

\proclaim{Lemma} 
Suppose $\X$ and $\Y$ are two \pendos{\rm ,} $x\in X$,  $y \in Y$ are two points{\rm ,} and
 $f:X \to R$ and $g:Y \to U$.
The one\/{\rm -}\/sided couplings
of the form $\hat\mu_{({{\cal T}}_x,{{\cal T}}_y,A)}$\/{\rm ,}\/ $A\in {\cal A}$\/{\rm ,}\/ are the 
extreme points of the
one\/{\rm -}\/sided couplings of $R_i({{\cal T}}_x)$ and $U_i({{\cal T}}_y)$ 
and any one\/{\rm -}\/sided coupling $\hat\mu$ of
$R_i$  and $U_i${\rm ,} $i<0${\rm ,}
is  of the form
 $$\int \hat\mu_{({{\cal T}}_x,{{\cal T}}_y,A)}\,dm_{(x,y)}(A)d\hat\mu(x,y)$$ for 
some family of
probability measures $m_{(x,y)}$ on ${\cal A}$.
\endproclaim

\demo{Proof} We only need to show that any one-sided coupling $\overline\mu$ of 
variables of the 
  form $R_i(h)$ and $U_i(h')$
 is of the form 
$$\overline\mu=\int \hat\mu_{(h,h',A)}\,dm(A)$$ for some measure $m$ on ${\cal A}$.
The proof is by induction.  
We first show this for a single variable.  This is equivalent to showing that any 
self-coupling
of uniform measure on $\K$ is an average of measures supported on graphs of
permutations.  To see this suppose $\hat\mu_0$ is such a self-coupling of $\K$.  The 
knowing
relation on $\K\times \K$ given by  $\hat\mu_0(i,j)>0$ satisfies the hypotheses of the 
Hall marriage lemma and hence
there is a bijective subrelation; i.e., $\hat\mu_0=\alpha\hat\mu_{\pi}+(1-\alpha)\hat\mu_1$
with $\alpha>0$ and $\mu_\pi$ supported on the graph of the permutation $\pi$. 
Repeating this for
$\mu_1$ and so on, we obtain a representation of the measure as an integral of 
measures supported
on graphs of permutations.   
Using the one-sidedness of $\overline\mu$ we complete the result inductively as
we conclude that  $\hat\mu$ is written uniquely as an integral of couplings which node by 
node sit on the graphs of
 permutations applied at each
node, i.e.\ sit on graphs of  tree automorphisms.
\enddemo 

We are ready to show that $\hat t$ and $\overline t$ agree.

\demo{Proof of Lemma~{\rm \ref{t-hatequalst-bar}}}
We begin   once more noting that by translating the random variables
$R_i$ and $U_i$ by $n$, both $\hat t_n$ and $\bar t_n$ are calculated as infima
over couplings of variables indexed on $-n\leq i<0.$  
 We have already noted that a pairing of the nodes of ${{\cal T}}_n$ by a
tree automorphism $A$ when viewed on the names along the branches of the
tree gives a one-sided coupling $\hat\mu_{(h,h',A)}$ of the distributions of names.  
This is enough to conclude that
$\hat t_n\leq \bar t_n$. 
As a one-sided coupling of $R_i({{\cal T}}_x)$ and $U_i({{\cal T}}_y)$ 
can be written as an integral over $X\times Y\times {\cal A}$
of couplings $\hat\mu_{({\cal T}_x,{\cal T}_y,A)}$ we see the other
inequality $\bar t_n\leq \hat t_n.$ 
\enddemo

\section{Tree finitely determined}
\label{tfd}

Now that we have discussed the theory of one-sided couplings we are
ready to define tree finitely determined.  This will play a major role
in the proof of the copying lemma.

 \numbereddemo{Definition} 
   We say  $(\X,f)$,  where $\X$ is a \pendo \
and $f:X\to R$ is tree adapted,
is {\it tree finitely determined}
(tree f.d.) if for every $\varepsilon>0$, there is a $\delta$ such that
 for any endomorphism $(Y,S,\mu)$ with function $g:Y \to R$
with 
 ${\rm tdist}(f,g)<\delta$ then 
 $\overline t(f,g)<\varepsilon.$ 
\enddemo 

 If $f$ is strongly tree adapted then tdist here can be replaced with dist.
This will enable us to  work with strongly
tree adapted functions.  Now we need the following lemma.  Later, as
we will see, tree f.d.\ and tree v.w.b.\ are equivalent.  We will see that all
tree-adapted factors of a tree f.d.\ map are tree f.d. At this point we need
something substantially less.

\proclaim{Lemma} 
  Suppose $((X,T,\mu),f)$ is tree {\rm f.d.}\ where $f$ is a generator.  Then for any
bounded map $h:X\to {\Bbb R}$, $((X,T,\mu), f\vee h)$ is tree {\rm f.d.}
\endproclaim

 \demo{Proof} 
This argument follows well established lines by approximation of $h$ by a
``finite code''.   Each successive step simply requires a closer match
in tdist.
  To begin, as $f$ is a generator, $h$ can be approximated arbitrarily well
in $L^1(\mu)$ by maps of the form 
$$H(f(x), f(T(x)),\dots,f(T^N(x)))$$ where $H$
is a continuous map from $R^{N+1}\to {\Bbb R}$ and $N$ is finite.  If $((Y,S,\nu),\break g\vee h')$
is sufficiently close in tdist to $((X,T,\mu),f\vee h)$ then $H(g(y),\break g(S(y)),\ldots, g(S^N(y)))$
will of necessity also be a good $L^1(\nu)$ approximation of $h'$.  

Now if $((Y,S,\nu),g\vee h')$ is close in tdist to $((X,T,\mu),f\vee h)$ then $((Y,S,\nu),g)$ is
close in tdist to $((X,T,\mu),f)$ but is not necessarily tree-adapted.  As $f$ is tree-adapted
though $g$ must separate inverse images of most points,
so some perturbation $g'$ of $g$ which agrees with $g$ on most of $Y$ will be tree-adapted.
If we are close enough in tdist then we will still have $H(g'(y),g'(S(y)),\dots,g'(S^N(y)))$\break a good 
approximation to $h'$ in $L^1(\nu)$ and $((Y,S,\nu),g')$ close in tdist to\break $((X,T,\mu),f)$, hence
close in $\overline t$.  Now if $((Y,S,\nu),g')$ and $((X,T,\nu),f)$ are close enough in $\overline t$ 
(notice that how close can be set after the value $N$ and continuous map $H$ are fixed) then
$((Y,S,\nu),g'\vee H(g',g'\circ S,\dots,g'\circ S^N))$ will be close in tdist to both
$((Y,S,\nu), g\vee h')$ and to $((X,T,\mu),f\vee h)$ which is to say $((X,T,\mu),f\vee h)$ is tree f.d.  
\enddemo 

\proclaim{{C}orollary} 
  Any tree {\rm f.d.}\ process has a tree {\rm f.d.}\ generator that is strongly tree adapted.
\endproclaim

\demo{Proof}   Choose $h$ in Lemma 4.2 to be a map to $\{1,2,\dots,p\}$ that
  separates inverse images.
\enddemo 

 \proclaim{Lemma} \label{treevwbimpliestreefd}
   If $(\X,f)$ is tree {\rm v.w.B.}\ then it is tree {\rm f.d.}
 \endproclaim 

\demo{Proof} 
Suppose $(\X,f)$ is tree v.w.B. Given $\varepsilon>0$ choose an $n$ so
that there exists a subset $X_0\subseteq X$ of measure $\geq 1-\varepsilon$ and
a fixed ${{\cal T}}_n,R$ name $h_{n}$
so that $\bar t_n({{\cal T}}_{x},h_{n})< \varepsilon$
for any $x\in X_0$.
For each $x \in X$ let $A_x$ be an automorphism that realizes the
minimum in the definition of $\bar t_n({{\cal T}}_{x},h_{n})$.
Using Lemma~\ref{dististdist}, choose a $\delta>0$ so small that if 
$\Y$ is a uniform\break $p$-to-one endomorphism, $g$ is a function to $R$, and
${\rm dist}(f,g)<\delta$
then there exists a subset $Y_0 \subset Y$ of measure greater than or equal to 
$1- 2\varepsilon$  such that  
$\bar t_n({{\cal T}}_{y},h_{n})<2\varepsilon$ 
 for all $y\in Y_0$.

Consider a  ${{\cal T}},R$
name $h$ constructed by tiling ${{\cal T}}$
with copies of $h_{n}$.  More precisely, 
for any $v$ such that $|v|=jn$ for some $j$ and any $v' \in \eta_n$
let $h(vv')=h(v')$.
For each $x \in X$ we will inductively construct an automorphism $A$
which will show ${{\cal T}}_x$ and $h$ are close in
$\bar t$.  The matching we use is a greedy algorithm matching $n$
levels at a time.
For each $v \in \eta_n$ let $A(v)=A_x(v)$.  Now assume $A$ has been
defined on all $v \in \eta_{jn}$.
For each $v \in \eta_{jn}$ and $v' \in \eta_n$ let 
$A(vv')=A(v)A_{T_{A(v)}(x)}(v')$.

Now we calculate $\bar t_{kn}({{\cal T}}_x, h)$.
Let $G_{kn}(x)$ be those nodes with $|v|=jn$ for some $0 \leq j <k$ 
and $T_v(x) \in X_0$.
Let $M(G_{kn}(x))$ be  the sum of $p^{-|v|}$ over all $v \in G_{kn}(x)$.
This construction leads to the calculation: 

$$\overline t_{kn}({{\cal T}}_x, h)\leq 1
	-M(G_{kn}(x))/k+ \varepsilon M(G_{kn}(x))/k.$$
The fact that 
$T$ is measure-preserving implies
 $$\int M(G_{kn}(x))\,d\mu(x)=k\mu(X_0) \geq k(1- \varepsilon).$$
Hence for all but $\sqrt{2\varepsilon}$ of the $x\in X$,
 $$\overline t_{kn}({{\cal T}}_x, h)
	\leq \sqrt{2\varepsilon}.$$
Precisely the same argument applied to $Y$ yields that for all but $2\sqrt{\varepsilon}$
of the $y\in Y$,
 $$\overline t_{kn}({{\cal T}}_y, h)\leq 2\sqrt{\varepsilon}.$$
We conclude that
 $$\overline t(f,g)=\liminf_{n\to\infty}\overline
t_n(f,g)\leq 4\sqrt{\varepsilon}$$ which ends the proof.
 \enddemo 
 
\proclaim{{C}orollary} 
  The standard endomorphism with the usual $p$ set 
independent generating partition is tree {\rm v.w.B}\ and hence tree {\rm f.d.}
\endproclaim
\demo{Proof} 
  It is trivial that with this partition every point $b$ has the
same past tree name.
Hence it
is tree v.w.B.
\enddemo 

We postpone the converse of Lemma \ref{treevwbimpliestreefd} 
as we will use the one-sided conjugacy 
theorem to prove it.

\section{The copying lemma}
\label{cl}

In this section we prove the isomorphism theorem for uniform $p$-to-one endomorphisms.
First we prove a Rokhlin lemma, then a strong Rokhlin lemma.  Next, we
prove the copying lemma.  The isomorphism theorem
will follow easily from the copying lemma.
Notice that the strong Rokhlin lemma is proved only 
for finite valued $f:X \to R$; i.e., $f$ is
a partition.  To emphasize this we refer to partitions $P$ and $Q$
instead of functions $f$ and $g$.

The first step is to prove a Rokhlin lemma for uniform endomorphisms.  
This result has appeared previously in the work of Rosenthal~\cite{Ro}.   We
present a proof here as his is perhaps too brief.

\numbereddemo{Definition} 
Let $\X$ be a uniformly $p$-to-one endomorphism 
and $T_v$ be some choice for the partial inverses of $T$.
A ${{\cal T}}_n$ {\it Rokhlin tree} is
a collection of disjoint sets $B_v\subseteq X$, $v\in\eta_n$, 
with the property that 
$T_{v}(B_\emptyset)=B_{v}.$  
\enddemo

\proclaim{Theorem}\label{Rokhlinlemma}
Let $\X$ be a uniformly $p$\/{\rm -}\/to\/{\rm -}\/one and ergodic endomorphism 
and $T_v$ be some choice for the partial inverses of $T$.
For each $n>0$ and $\varepsilon>0$ there exists a ${{\cal T}}_n$
Rokhlin tree $B_v$ with $\mu(\cup_v B_v)>1-\varepsilon.$
\endproclaim 

\demo{Proof}  
For any set $C$ and $n>0$ define 
$$B_{\emptyset} = \{ x \ | \ \min \{i \geq 0 : T^i(x) \in C\}= 0 
	\bmod (n+1) \} \setminus (\cup_{0 \leq i \leq n} T^{i}(C)).$$
Suppose $x \in B_{\emptyset} \cap T_v(B_{\emptyset})$ for some 
$v\in \eta_n$.  Then since 
$\mathbold{\cup}_{0 \leq i \leq n}T^{-i}(x) \cap C = \emptyset$
for any point $x \in B_{\emptyset}$,
$$\min \{i\geq 0 : T^i(T_v(x)) \in C\}-\min \{i\geq 0 : T^i(x) \in C\}=|v|.$$
Both of the terms on the left-hand side cannot be equal to $0$ mod
$(n+1)$ unless $|v|=0$.
Thus
$B_{\emptyset} \cap T_v(B_{\emptyset}) =\emptyset$
for any $0 < |v| \leq n$ which implies
$B_v=T_v(B_\emptyset)$ forms a Rokhlin tree.  

Since 
$\mu(\{ x \ | \ (\min_{i\geq 0} T^i(x) \in C)= i\})$ is nonincreasing  
we have
$$\mu(B_\emptyset)>1/(n+1) - (n+1)p^{n+1}\mu(C)\hbox{ and }\mu(\cup_v B_v)>1 
	- (n+1)^2p^{n+1}\mu(C).$$
This last term 
can be made as small as we like by choosing $\mu(C)$ small.  
\enddemo 

Now we prove a strong Rokhlin lemma.  This says that the top
level $B_\emptyset$ can be chosen independently of any partition.

\proclaim{Lemma} 
\label{strongrokhlin}
Suppose $(X,T,\mu)$ is an ergodic uniformly $p$\/{\rm -}\/to\/{\rm -}\/one
endomorphisms and $P$ is a finite partition of $X$.  For any
$\varepsilon>0$ and $n$ there is a ${{\cal T}}_n$\/{\rm -}\/Rokhlin tree $C_{v}$
so that $$ \mu\left( \cup_{v \in \eta_n}C_v\right) >1-\varepsilon$$ and
$${\rm Dist}(P)={\rm Dist}(P|C_{\emptyset});$$ i.e.{\rm ,} $P$ and $C_{\emptyset}$ are independent.
\endproclaim 

\demo{Proof}  
Given $\varepsilon$ choose $m>4(n+1)/\varepsilon.$
Let $B_{\emptyset}$ be the top of a ${{\cal T}}_m$ Rokhlin tree with
$\mathbold{\cup}_{v \in \eta_m}B_v >1-\varepsilon/2$.
For each $$A \in (\vee_{v \in \eta_m} T_v(P))|_{B_{\emptyset}}$$
write $A= \cup_{0}^{n}A_i$ with each $\mu(A_i)=\mu(A)/(n+1)$.
Now let 
$$C_{\emptyset}=
  \mathbold{\cup} _{A \in (\vee_{v \in \eta_m} T_v(P))|_{B_{\emptyset}}} 
  \mathbold{\cup} _{j=0}^{m/(n+1)-2} \mathbold{\cup} _{i} T^{-j(n+1)-i}(A_i).$$
Then $C_{\emptyset}$ forms the top of a ${{\cal T}}_n$ Rokhlin tree with
$\mathbold{\cup} _{v \in \eta_n}C_v >1-\varepsilon$.
We also have
$$\Prob(P\ | \ \mathbold{\cup} _{i=0}^{m-n-1} T^{-i}(B_\emptyset))=\Prob(P\ | \ C_\emptyset).$$
Thus for any element $P_i$ of $P$
$$\mu(P_i \cap C_\emptyset)>(1-\varepsilon)\mu(P_i)/(n+1). $$
Then we pare down $C_\emptyset$ to $C'_\emptyset$ so that for every 
element $P_i$ of $P$ we have precisely
$$\mu(P_i \cap C'_\emptyset)=(1-\varepsilon)\mu(P_i)/(n+1) $$
and we are done.
\enddemo

Now we are ready to prove the copying lemma which is the main element
in the proof of Theorem~5.5.

\proclaimtitle{copying lemma}
\proclaim{Lemma} \label{copyinglemma}
Suppose $(X,T,\mu)$ and $(Y,S,\nu)$ are \pendos\ with finite
tree adapted functions $f$ and $g$ and $\hat \mu\in
J_e^+$ an ergodic joining.  For any $\varepsilon>0$ there is
a tree adapted function $\tilde f:Y \to R$ with
$${\rm dist}(f\coup{\hat \mu} g,{\tilde f}\coup{\nu} g)<\varepsilon.$$ 
As a consequence{\rm ,} if $R_i$ is tree {\rm f.d.}\ then 
 $${\cal O_{\varepsilon'}}=
	\{\tilde\mu\in J^+_e|f\cont{\varepsilon'}{\tilde\mu}Y\}$$
is an open and dense subset of $J^+_e.$
\endproclaim 

\demo{Proof} 
First we show the result for finite partitions $P$ and $Q$ instead of
functions $f$ and $g$.
The definition of dist gives   an $n$ and a
$\delta$ such that if 
$$\sum_{A \in \vee^{n}_{0}(T\times S)^{-i}(P \vee Q)}
	{|\nu(A)-\hat\mu(A)|}< 1 - \delta,$$
then 
$${\rm dist}(P\coup{\hat \mu} Q,{\tilde P}\coup{\nu} Q)<
	\varepsilon.$$ 
For any partition $P$ and any $m$ define a new partition $P^m$ so that $x$ and $x'$
are in the same element of $P^m$ if and only if 
$\bar t_{m}({{\cal T}}_x,{{\cal T}}_{x'})=0$.
Choose
$m > 4n/ \delta$. Build ${{\cal T}}_m$ 
Rokhlin trees $B_v \subset X$ and $C_v \subset Y$ so that
$\mu(\cup_{v}B_v)=\nu(\cup_{v}C_v) >1-\delta/2$, $B_\emptyset$ is
independent of $P^m$, and $C_\emptyset$ is independent of $Q^m$.

Next define a measure-preserving 
map $h:B_\emptyset \to C_\emptyset$ so that the measure
$h$ generates on $B_\emptyset \times C_\emptyset$ restricted to 
$P^m \vee Q^m$ is the same as the measure $\hat \mu$
restricted to $P^m \vee Q^m$.  This is possible because $B_\emptyset$ is
independent of $P^m$, and $C_\emptyset$ is independent of $Q^m$.

As $\hat\mu$ is one-sided we can write it as 
$$\hat\mu=\int\mu_{({\cal T}_x,{\cal T}_y,A)}\,dm_{(x,y)}d\hat\mu(x,y).$$
Lift $\hat\mu$ to $X\times Y\times {\cal A}$ as
$$\bar\mu=\int\delta_x\times\delta_y\times m_{(x,y)}\,d\hat\mu(x,y).$$
Now $P^m\times Q^m\times {\cal A}_m$ is a finite partition of
$X\times Y\times{\cal A}$.  As $C_\emptyset$ is independent of
$Q^m$ we can  extend $\nu$ on $C_\emptyset$ to $\bar\nu$
on $C_\emptyset\times{\cal A}$ to be identical in distribution 
on $Q^m\times{\cal A}_m$ to $\bar\mu$ (when normalized).
We can now construct a measure-preserving map $\bar h:(B_\emptyset,\mu)\to (C_\emptyset\times {\cal A},\bar\nu)$
so that the normalized measure supported on the graph of $\bar h$ 
restricted to $P^m\times Q^m\times 
{\cal A}_m$ agrees in distribution with $\bar\mu.$

Now we are ready to define the new partition $\tilde P$.  
Write $\bar h(x)=(h(x),A(x)).$
For each $x \in B_\emptyset$ and every $v \in {{\cal T}}_m$ set
$${\tilde P}(T_{A(x)(v)}h(x))=P(T_v(x)).$$
On the rest of the space define $\tilde P$ in any way such that 
$\tilde P$ is tree-adapted.

Now for any set $A \in \vee_{0}^{n} (T\times S)^{-i}(P \vee Q)$,
$$\nu(A \cap (\cup_{|v|<m-n}C_v))=\hat\mu(A \cap
	(\cup_{|v|<m-n}B_v)).$$
Since $\hat\mu(\cup_{|v|<m-n}B_v) = \nu
	(\cup_{|v|<m-n}C_v)>1-\delta$ this implies
$$\sum_{A \in \vee_{0}^{n} (T\times S)^{-i}(P \vee Q)}
		{|\nu(A)-\hat\mu(A)|}< 1 - \delta.$$
Thus
\vglue-4pt
\centerline{${\rm dist}(P\coup{\hat \mu} Q,{\tilde P}\coup{\nu} Q)<
	\varepsilon.$}  
\vglue6pt\noindent  This completes the proof of the first statement for finite partitions.

To extend the result to strongly tree-adapted functions do the following.
Partition $R$ into sets of small diameter.  Choose one representative point in
each element of the partition. Define $F:R\to R$ to map all points in each
partition element to its representative point.  It follows that if the partition
is fine enough then ${\rm dist}(f,F\circ f)$ will be small.  Do the same with
$g$, i.e.\ construct a finite valued $G:U\to U$.  Once more if the partitions are
fine enough then we will have 
$${\rm dist}(f\coup{\hat\mu} g,F\circ f\coup{\hat\mu}G\circ g)<\varepsilon/3.$$
As $f$ and $g$ are strongly tree-adapted, if the partition elements are less than 1
in diameter $F\circ f$ and $G\circ g$ will also be strongly tree-adapted (by
the discrete metric on their finite ranges).   We can now apply the finite partition
version proven above to $F\circ f$ and $G\circ g$ to construct $\tilde f$ 
with  
$${\rm dist}(F\circ f\coup{\hat\mu}G\circ g,\tilde f\coup{\nu}G\circ g)<\varepsilon/3,$$
knowing $G$ comes from a fine enough partition, independent of the choice of $\tilde f$
that
$${\rm dist}(\tilde f\coup{\nu}G\circ g,\tilde f\coup{\nu}g)<\varepsilon/3.$$ 
 
For the second result notice that the previous equation implies
$${\rm dist}(f ,{\tilde f})<
	\varepsilon.$$ 
For any $\delta$
if $(\X,f)$ is tree f.d.\ then $\varepsilon$ and $\tilde f$ can be chosen so that 
$$\bar t (f,\tilde f)<\delta.$$ 
For any $\delta_1$ if $\delta$ and $\varepsilon'$ are small enough then
by Lemma~\ref{partialsgroup} we can extend $\nu$ to
a one-sided and ergodic joining $\tilde \mu$ with
$${\rm dist}(\tilde f\coup{\nu} g ,f\coup{\tilde\mu}g)<\delta_1\hbox{ and } f\cont{\delta}{\tilde\mu}Y$$
where 
${\cal O_{\varepsilon'}}$ is dense.  Openness follows easily from
the definition.
\enddemo

\proclaim{Theorem}\label{mainresult}
  A \pendo\ $\X$ is one\/{\rm -}\/sidedly conjugate to the standard \pendo\ 
$(B,\sigma,\nu)$ 
 if and only if there exists a generating function $f$ so that
$\X$ and $f$ are tree {\rm v.w.B. (}\/or equivalently tree {\rm f.d.).}
\endproclaim 

\demo{Proof} 
Let $f$ be a strongly tree adapted function from $X$ to $R$ 
and $K$, the
standard independent generating partition of $B$.
We know $J_e^+$ is a $G_\delta$ subset of $J^+$ in the weak* topology.
Since $\X$ and $f$ are tree v.w.B. they are also tree f.d.  
Thus 
the copying lemma tells us that the 
${\cal O}_\varepsilon$ are open and dense in $J_e^+$.
Intersecting over $\varepsilon=1/n$, we see that the Baire category theorem shows that the set of $\hat\mu$ with
$f\cont{0}{\hat\mu} B$ is a residual subset of $J_e^+$.  
Let $K$ be the standard independent generating partition of $B$.
As the standard example 
is also tree v.w.B. and tree f.d.,
the set of $\hat\mu$ with $K\cont{0}{\hat\mu} X$ is a residual
subset of $J_e^+$. Thus 
the set of $\hat\mu$ with $K\cont{0}{\hat\mu} X$
and $f\cont{0}{\hat\mu} B$ is nonempty and $(X,T,\mu)$ and $(B,\sigma,\nu)$ are
isomorphic.

Notice  that we need only assume $f$  tree adapted here and not necessarily strongly, 
as we can extend $f$ to an $f\vee h$ which is still tree f.d. and now is
strongly tree adapted and hence isomorphic to the standard example.

   All that remains is to show tree f.d.\ and tree
v.w.B. are equivalent. We already know that tree v.w.B. implies tree f.d.
For the other implication, we have just seen that tree f.d.\ implies
one-sidedly conjugate to the standard action which is tree v.w.B.   We
have also seen that tree v.w.B. descends to one-sided and tree adapted
factors. Hence any \pendo\ which is 
isomorphic to a tree v.w.B. endomorphism is
tree v.w.B.   Thus tree f.d.\ implies tree v.w.B.   
\enddemo

\section{Examples of tree v.w.B. skew products}\label{tid}

We will show now that   general classes of isometric extensions
of standard endomorphisms are all one-sidedly Bernoulli.
Among these will be the $\tid$ endomorphisms, where $T$ is
an irrational rotation. These were described in the
first section.  

Throughout this section we consider $p$ to be
fixed.  Remember that $(B,\sigma,\nu)$ is 
the standard \pendo.
We also fix  $(Z,d)$  a compact metric space with ${\cal  I}$ 
its space of isometries.  We assume ${\cal  I}$ acts transitively,
i.e.\ $Z$ is a homogeneous space, and has  Haar measure. 
We put on ${\cal  I}$ the uniform
topology.  
Given a function $f:B\to {\cal  I}$ we construct
the {\it cocycle extension} $T_f$ acting on $B \times Z$
by
 $$T_f(b,z)=(\sigma (b),f(b)(z)).$$

The map $\sigma$ has a natural set of partial inverses $\sigma_v$.
These extend to form a natural set of partial inverses $(T_f)_v$.
Set $c_m=\sup_{|v|=m } ({\rm diam}(f(\sigma_v(B)))).$ 

\numbereddemo{Definition} 
We say $f$ {\it generates a summable cocycle} if $\sum_m c_m<\infty.$
\enddemo 

Notice that any $f$   depending on only finitely many coordinates $b_0,\dots,b_m$
automatically generates a summable cocycle.  
Thus the $\tid$ endomorphisms generate summable cocycles.

\proclaim{Lemma} 
\label{lem1}
Suppose $f:B\to {\cal  I}$ generates a summable cocycle.  Then for
any $\varepsilon>0$ there exists $\delta$ such that if 
$d((b,z),(b',z')) < \delta$ then $d((T_f)_v(b,z),\break (T_f)_v(b',z')) < \varepsilon$
for all $v \in \eta$.
\endproclaim 

\demo{Proof} 
Given $\varepsilon$ there exists $n_0$ so that 
$\sum_{n\geq n_0}{c_n} \leq \varepsilon/2$.  Choose 
$\delta < \varepsilon/2$ such that $d(b,b')<\delta$ implies that 
$b_i=b'_i$ for all $0 \leq i <n_0$. 
By the definition of $c$ if $b_i=b'_i$ for all $0 \leq i <n_0$ then
\vglue9pt
\hfill ${\displaystyle d((T_f)_v(b,z),(T_f)_v(b',z')) \leq d((b,z),(b',z')) + \sum_{n\geq n_0}{c_n}
	\leq \varepsilon .}$
\enddemo 

\proclaim{Lemma} 
\label{lem2}
Suppose $f:B\to {\cal  I}$ generates a summable cocycle and $T_f$
is weakly mixing.  For any $\varepsilon>0$ 
there is an $N>0$ so that
$(T_f)^{-N}(b,z)$ is $\varepsilon$ dense in $B\times Z$ for all $(b,z)$ .
\endproclaim 

\demo{Proof} 
From \cite{R} we know that if $T_f$ is weakly mixing then it must be
v.w.B. and hence $K$.
The fact that
$T_f$ is a $K$-system implies that there exists an $N'$ such that 
for most $b$ and all $z$ we have that 
$(T_f)^{-N'}(b,z)$ is $\varepsilon/2$ dense in $B\times Z$.
In particular this holds for one $b$.
Lemma \ref{lem1} implies that there exists $N''$ such that 
if $b_i=b'_i$ for all $0 \leq i <N''$ then
$(T_f)^{-N'}(b',z)$ is $\varepsilon$ dense in $B\times Z$.
For every $b''$ there exists $\hat b \in \sigma^{-N''}b''$ such that 
$\hat b_i=b_i$ for all $0 \leq i <N''.$
Thus 
$(T_f)^{-N'-N''}(b',z)$ is $\varepsilon$ dense in $B\times Z$.
\enddemo 

\proclaim{Theorem}
Suppose that $f$ generates a summable cocycle
and that $T_f$ is weakly mixing.  Then $T_f$ is tree {\rm v.w.B. }
\endproclaim 

\demo{Proof} 
Given $\varepsilon$ we get a $\delta$ from Lemma \ref{lem1} and an $N$
from Lemma \ref{lem2} which implies that 
$(T_f)^{-N}(b,z)$ is $\delta$ dense in $B\times Z$ for all
$(b,z) \in B \times Z$.
For any $(b,z),(b',z') \in B \times Z$
we define the tree automorphism $A$ that pairs ${{\cal T}}_{(b,z)}$ with
${{\cal T}}_{(b',z')}$ inductively. We pair at least $N$ levels at a time.  
If $d((b,z),(b',z'))<\delta$ then by Lemma
\ref{lem1} we are done.  Otherwise Lemma \ref{lem2} implies that there
exists 
a tree automorphism ${{\cal A}}_N$ which pairs at least one preimage of
$(b,z)$ with a preimage of $(b',z')$ so that the pair is within $\delta.$

Now suppose we have defined $A$ up to at least level $kN$.  We will extend it
to at least level $(k+1)N$.  
If $A(vv')$ has been defined for all $v' \in \eta $ then we need
to do nothing.
If it has not and
$d((T_f)_v(b,z), (T_f)_{A(v)}(b',z'))<\delta$  
then extend it by the identity automorphism.  
By this we mean, for all $v' \in \eta$, that $A(vv')=A(v)v'.$
If neither of the above conditions is satisfied then 
we use Lemma \ref{lem2} to tell us how to define 
$A(vv')$ for all $v' \in \eta_N $.

Choose $k$ so that $(1-1/p^N)^{\varepsilon k} < \varepsilon.$
By the previous paragraph and Lemma~\ref{lem1} we know that
for each $n>\varepsilon kN$ the fraction of preimages that is paired
within $\varepsilon$ is at least $1-\varepsilon$. Thus
$${\bar t}_{kN} ((b,z),(b',z')) < 3 \varepsilon$$
and $T_f$ is tree v.w.B. 
\enddemo 

As well as the $\tid$ examples, our methods cover the following smooth
endomorphisms. 
Replace $B$ by $x\to 2x$ on ${\Bbb R}/{\Bbb Z}$ and let $Z={\Bbb R}/{\Bbb Z}$ as well.  Set $f(x)$ to be any H\"older function to ${\Bbb R}/{\Bbb Z}$
that is not a coboundary (for example $f(x)=\sin(2\pi x)$) giving a smooth action
on the 2-torus.  That $f$ is not a coboundary means $T_f$ is a weakly mixing action
and that $f$ is H\"older implies that it generates a summable cocycle.
Thus we conclude that such an action must be tree
v.w.B.   It would be interesting to know if there can be a smooth and
uniformly $p$-adic action that is v.w.B. but not tree v.w.B.   Our work
here shows it will not be found among the isometric extensions of
$x\to px\, {\rm mod}\,  1$ with summable cocycles.

\end{document}